\documentclass[a4paper,10pt]{article}
\usepackage{amsmath, amsthm, amssymb}
\usepackage{graphicx}

\title{Nonalternating knots and Jones polynomials}
\author{Neil R. Nicholson\\
	The University of Iowa\\
	\\
	http://www.math.uiowa.edu/$\sim$nrnichol\\
	}

\newtheoremstyle{dotless}{}{}{\itshape}{}{\bfseries}{}{ }{}
\def\Real{\hbox{I\kern-.1667em\hbox{R}}}
\theoremstyle{dotless}	
\input epsf

\newtheorem{theorem}{Theorem}[section]
\newtheorem{corollary}[theorem]{Corollary}
\newtheorem{lemma}[theorem]{Lemma}
\newtheorem{conjecture}[theorem]{Conjecture}
  
\newcommand{\kb}[1]{\langle #1 \rangle}

\begin{document}

\maketitle

\begin{abstract}
We consider here nonalternating knots and their properties.  Specifically, we show certain classes of knots have nontrivial Jones polynomials.
\end{abstract}

\section{Introduction}

It is known that there are nontrivial links with unit Jones polynomial \cite{15} as well as a nontrivial virtual knot with trivial Jones polynomial for each nontrivial standard knot \cite{5}.  Moreover, there have been numerous pairs of distinct knots found to have the same Jones polynomial \cite{17}.  Yet it is still unknown if the Jones polynomial detects unknottedness.

Bae and Morton \cite{2} developed a simple combinatorial method for calculating the potential two extreme terms of the Kauffman bracket polynomial and consequently of the Jones polynomial for unoriented diagrams.  In Section \ref{jones} of this paper we will use their tool to place a stronger bound on the span of the Jones polynomial in terms of the nonalternatingness of a diagram.  In Section \ref{classes} we prove that numerous classes of knots have nontrivial Jones polynomials, including untwisted Whitehead doubles of any knot.  Our proof is different than that given by Stoimenow \cite{13} for untwisted Whitehead doubles of positive knots.  Other knots we show to have nontrivial Jones polynomials include any Whitehead double or cable knot of an alternating knot, all pretzel knots, and any knot of $o$-length two (and consequently any almost alternating knot).

\section{Preliminaries}
\label{preliminaries}

Throughout this paper we will assume that $K$ is a knot (unoriented, unless otherwise noted) in general position with diagram $D$ and universe (shadow) $U$.  Travelling along $D$, if we encounter $m+1$ consecutive overcrossings, we will call that portion of the diagram an \textit{overpass of length m} (sometimes referred to as a \textit{bridge of length m}).  An edge of $U$ ($U$ is a 4-valent graph) is said to be \textit{positive} if in $D$ it is part of an overpass and is bounded by two overcrossings.  The number of positive edges in $D$ will be defined as the \textit{o-length} $o(D)$ of $D$.  Define $o(K)$ to be the minimum $o$-length over all embeddings of $K$.  \textit{Negative edges}, $u$-length, $u(D)$, and $u(K)$ are defined similarly using undercrossings.  The bridge number $b(D)$ is the number of overpasses of nonnegative length, and the bridge index $b(K)$ is the minimum such value over all diagrams of $K$.

\begin{lemma}
\cite{14} The boundary of any region in the complement of $U$ has an equal number of positive and negative edges. \label{equalposneg}
\end{lemma}

\begin{corollary}
$o(D)$ = $u(D)$. \label{oequalsu}
\end{corollary}

\begin{lemma}
$o(K_{1} \# K_{2}) \leq o(K_{1}) + o(K_{2})$. \label{olengthconnectedsum}  
\end{lemma}

Notice if both the $K_{i}$ and the connected sum in the above lemma are alternating than the equality holds.

\begin{lemma}
$b(D) + o(D) = c(D)$. \label{bridgeolength}
\end{lemma}

Lemma \ref{bridgeolength} does not generalize to indices.  In fact, the equality fails for all nontrivial $2$-bridge knots $K$, since $b(K) = 2$, $o(K) = 0$ ($2$-bridge knots are alternating), and $c(K) \geq 3$.\\

It is easily seen that the Gauss code for a diagram and the $o$-length are related.

\begin{lemma}
$o(D)$ equals the number of consecutive positive pairs (considered cyclically) in the Gauss code for $D$.
\end{lemma}

Recall that a rational link $K = C(a_{1}, a_{2}, ..., a_{n})$ is formed by taking two parallel, horizontal strings (if $n$ is odd, take the strings to be vertical) and performing $a_{1}$ half-twists to the right (resp. bottom) ends of the strands.  Now perform $a_{2}$ half-twists to the bottom (resp. right) ends of the strands.  Continue alternating between the right and bottom ends of the strings.  Upon concluding, identify the upper two ends as well as the lower two ends.  These diagrams appear to have nonalternation, specifically occurring between twists of opposite sign.  However, if $K$ is rational (considering $K$ to be of one component), we know $o(K) = 0$, as each rational knot is equivalent to an alternating and hence adequate knot \cite{6}.  Moreover, as the set of rational knots is equal to the set of knots of bridge index two, Lemma \ref{bridgeolength} does not generalize to indices for any rational knot.

Notice that $o(K)$ gives us a sense of the alternatingness of a given knot type.  Knot tables order knots based on their crossing index and within each index by whether or not the knot is alternating.  It seems only natural to proceed in grouping the nonalternating knots based on ``how nonalternating" they really are; that is, by $o(K)$.  It is interesting to note that many of the nonalternating $8$- and $9$-crossing knots in Rolfsen's table \cite{11} are illustrated using diagrams of $o$-length three when each actually has $o$-index two (as we note later, no knot can have $o$-length one).  Of the nonalternating knots of crossing index less than ten, all but one, the Perko knot, have been found to have $o$-index two.

\begin{conjecture}
$o(10_{161}) = 3$
\end{conjecture}

Note that $10_{161}$ is also the first knot from Rolfsen's table which strictly satisfies an inequality proven by Kidwell that relates the maximum degree of the $Q$ polynomial and maximal overpass length \cite{12}.

There are two tools that will help us prove our main theorems.  The first and more widely known of these is the Kauffman bracket polynomial.  For its original definition and properties, see \cite{4}.  We give a brief review of its definition here.

Any crossing of a diagram locally separates the plane into regions, as follows.  If we rotate the overcrossing strand counterclockwise, the regions swept out will be referred to as the \textit{A} regions.  Rotating the strand clockwise sweeps out the $B$ regions.  An $A$-smoothing of the crossing results by replacing the crossing with two smooth edges so that the $A$ regions are connected.  We also have $B$-smoothings.

A \textit{state} of $D$ is a choice of smoothing for each crossing.  Let $a(S)$ and $b(S)$ be the number of $A$- and $B$-smoothings in the state $S$, respectively. Define $S_{A}$ (resp. $S_{B}$) to be the state resulting from replacing all crossings with $A$-smoothings (resp. $B$-smoothings).

For a diagram $D$, the \textit{Kauffman bracket polynomial} is a Laurent polynomial in the variable $A$ with integer coefficients.  It is given by: 

\begin{center}
$\kb{D} = \sum A^{a(S)-b(S)} (-A^{2} - A^{-2})^{|S|-1}$,
\end{center} 

\noindent where the sum is taken over all states $S$ and $|S|$ is the number of simple closed curves gotten by smoothing each crossing according to $S$.  The Kauffman bracket is an invariant of regular isotopy (invariant under Reidemeister II and III moves).  If $D$ is oriented, assign a value of $+1$ or $-1$ to each crossing according to the usual right-hand rule.  Define the writhe of $D$ $w(D)$ to be the sum of these values.  Multiplication by a factor of $(-A)^{-3w(D)}$ and substituting $A = t^{-1/4}$ yields the Jones polynomial $V_{K}(t)$, an ambient isotopy invariant (invariance under all Reidemeister moves).

Let $S$ be any state.  We then define $max(S) = a(S) - b(S) + 2|S| - 2$ and $min(S) =  a(S) - b(S) - 2|S| + 2$.

\begin{lemma}
\cite{2} The maximum degree of $\kb{D}$ is less than or equal to \linebreak$max(S_{A}) = c(D) + 2|S_{A}| -2$.  Similarly, the minimal degree is greater than or equal to $max(S_{B}) = -(c(D) + 2|S_{B}| -2)$. 
\end{lemma}

Thus it makes sense to define $a_{S_{A}}$ and $b_{S_{B}}$ to be the coefficients (potentially zero) of the terms of  degree $c(D) + 2|S_{A}| -2$ and -($c(D) + 2|S_{B}| -2$) in $\kb{D}$, respectively.

\begin{figure}[hbtp]
  \begin{center}
    \leavevmode
    \epsfxsize = 4.75cm
    \epsfysize = 3cm
    \epsfbox{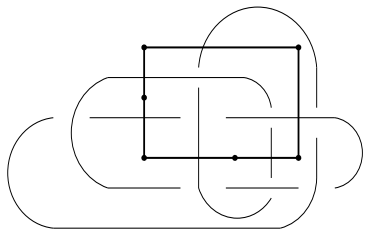}
    \caption{$8_{20}$ with $G$ superimposed}
    \label{knot820withG}
  \end{center}
\end{figure}

The second tool we utilize comes from a paper by Bae and Morton \cite{2}.  The complement of $U$ consists of disjoint regions of the plane.  Place a vertex of $G$ in each region whose boundary contains a positive edge (or equivalently by Lemma \ref{equalposneg}, a negative edge).  If $e$ is a positive or negative edge of $D$, it on the boundary of exactly two regions in which we have placed vertices.  Connect these vertices by an edge $e'$ of $G$ so that if $G$ is superimposed on $D$, $e'$ intersects $e$ transversally.  Each edge of $G$ naturally inherits a sign from $D$.  Figure \ref{knot820withG} shows a diagram of $8_{20}$ with $G$ superimposed.  $G$ is called the \textit{nonalternating skeleton of $D$} and need not be connected.  

Following immediately from Lemma \ref{oequalsu} and the construction of $G$ we have:

\begin{lemma}
The nonalternating skeleton for any knot diagram $D$ has 2$o(D)$ edges.
\end{lemma}

\section{Jones polynomials}
\label{jones}

For nonalternating knots, span($V_{K}(t)$) $\leq$ $c(D)$ - 1 for any totally reduced diagram (defined below) $D$ of $K$ \cite{16}.  Using the nonalternating skeleton, however, we are able to strengthen this bound.  We first describe more of the construction and results from Bae and Morton.  The symbols $K$, $D$, $U$, and $G$ will be as in the previous sections.

Each vertex of $G$ has even valency.  Moreover, the edges of $G$ must alternate sign around each vertex.  Therefore we can locally replace each $2m$-valent vertex of $G$ by the $m$ $2$-valent vertices formed by connecting each positive edge with its counterclockwise (resp. clockwise) neighboring negative edge.  The resulting collection of simple closed curves will be referred to as $G_{A}$ (resp. $G_{B}$).

Suppose we have $D$ with $G_{A}$ superimposed.  Crossings of $D$ locally separate the plane into $A$ and $B$ regions, as previously described.  Call any arc that intersects $D$ at a single crossing $c$, approaching $c$ through the $B$ regions formed by $c$, and whose endpoints lie on a single component of $G_{A}$ a \textit{Lando-$B$ arc}.  Form $L_{A}$, \textit{Lando's $A$-graph}, by taking a vertex for each Lando-$B$ arc and an edge connecting two vertices if and only if the endpoints of the corresponding arcs alternate around the same component of $G_{A}$.  Similarly, form $L_{B}$ using $G_{B}$ and Lando-$A$ arcs.  $L_{A}$ is used when calculating $a_{S_{A}}$ as $L_{B}$ is used to calculate $b_{S_{B}}$.

\begin{figure}[hbtp]
  \begin{center}
    \leavevmode
    \epsfxsize = 2.5cm
    \epsfysize = 2.5cm
    \epsfbox{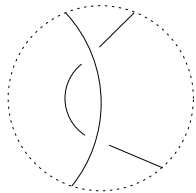}
    \caption{A II-reducible portion of a knot}
    \label{reducible}
  \end{center}
\end{figure}

$D$ is called \textit{reduced} if there is no circle in the plane meeting the corresponding shadow transversally at a crossing and not intersecting the diagram in any other place.  $D$ is called \textit{II-reduced} if there are no obvious removable Reidemeister-$II$ moves; i.e., the knot contains no $2$-tangle as in Fig. \ref{reducible}.  If $D$ is both reduced and $II$-reduced, we will call $D$ \textit{totally reduced}.

Recall that a \textit{simple graph} is a graph which contains no loops or multiple edges.

\begin{lemma}
The nonalternating skeleton of any prime knot can be considered to be a simple graph.  \label{nosharedvert}
\end{lemma}
\begin{proof}
Obviously there are no loops.  Suppose two distinct edges $e_{1}$ and $e_{2}$ of $G$ share endpoints $v_{1}$ and $v_{2}$.  Then $e_{1}e_{2}$ forms a simple closed curve and hence separates the plane into two regions.  We can consider each region as a $1$-tangle, since only one edge of $D$ intersects each of $e_1$ and $e_2$.  The closure of each of these tangles forms a knot.  If both $1$-tangles form nontrivial knots when closed, then we contradict $K$ being prime.  Thus, one of the tangles, when closed, forms a trivial knot.  Unknot this portion of the knot and the result follows.
\end{proof}

As a direct consequence of the above lemma, we have the following:

\begin{theorem}
If $D$ is a totally reduced diagram for a prime knot $K$, then \linebreak$o(D) = 0$ if $K$ is alternating and we can assume $o(D) \geq$ 2 if $K$ is nonalternating.
\end{theorem}

Let $|G_{A}|$ and $|G_{B}|$ be the number of such curves for these collections.  Then \cite{2} proves the following \textit{extreme states bound} for span$(V_{K}(t))$.

\begin{theorem}
\cite{2} If $c(D) = n$ and $v(G)$ is the number of vertices in $G$, then 
span($V_{K}(t)$) $\leq$ n + $\frac{1}{2}$($\mid$$G_A$$\mid$ + $\mid$$G_B$$\mid$ - $v$($G$)). \label{extremestatesbound} 
\end{theorem}

This bound helps us prove one of our main results: a bound on the span of the Jones polynomial in terms of $o$-length.  We first prove a necessary lemma:

\begin{lemma}
If $D$ is a totally reduced diagram of odd $o$-length with $o(D) \geq 5$ and $G$ its nonalternating skeleton, then $|G_{A}| + |G_{B}| - v(G) \leq -6$.
\end{lemma}
\begin{proof}
For $o(D) = 5$, the result holds by simply inspecting the possible nonalternating skeletons.  Suppose for a general $m > 2$, if $o(D) = 2m + 1$, then $|G_{A}| + |G_{B}| - v(G) \leq -6$.  To obtain a skeleton $G'$ for a diagram $D'$ of $o$-length $2m + 3$, we must add an additional four edges to $G$.  Note that all such skeletons can be built up this way, allowing the induction to proceed.

In order to preserve the structure of $G$, edges must be added to $G$ in pairs as a single segment (one positive edge and one negative edge).  We must consider adding two such segments $S_{1}$ and $S_{2}$ to $G$, and there are three ways that this can be done:\\  

\noindent \textit{Method 1}: The first method involves replacing a $2$-valent vertex of $G$ by a segment.  If one of $S_{1}$ or $S_{2}$ is added this way, then the other must be as well, though not necessarily near the other in $G$.  Notice that no extra circles result when forming $G_{A}$ or $G_{B}$.  By the inductive hypothesis we have
\begin{align*}
|G'_{A}| + |G'_{B}| - v(G') &= |G_{A}| + |G_{B}| - (v(G) + 2) \\
&\leq -8
\end{align*}

\noindent \textit{Method 2}: The new edges can be added by placing the endpoints of $S_{1}$ on different vertices of $G$.  Again by Lemma \ref{nosharedvert}, $S_{2}$ must be added this way and its endpoints must lie on the same vertices as $S_{1}$.  Notice that $|G_{A}|$ and $|G_{B}|$ increase by at most one, giving:
\begin{align*}
|G'_{A}| + |G'_{B}| - v(G') &\leq |G_{A}| + |G_{B}| - (v(G) + 2) \\
&\leq -6
\end{align*}

\noindent \textit{Method 3}: The third method is obtained by forming a square from $S_{1}$ and $S_{2}$ disjoint from $G$.  Then by the inductive hypothesis
\begin{align*}
|G'_{A}| + |G'_{B}| - v(G') &= |G_{A}| + 1 + |G_{B}| + 1 - (v(G) + 4) \\
&< -8
\end{align*}

In all cases, we have $|G'_{A}| + |G'_{B}| - v(G') \leq -6$.
\end{proof}

\begin{theorem}
Let $D$ be a totally reduced $n$-crossing diagram for a prime knot $K$.  Then:
\begin{enumerate}
\item if $o$($D$) = 3, then span($V_{K}(t)$) $\leq$ n-2
\item if $o$($D$) is even, then span($V_{K}(t)$) $\leq$ n-1
\item if $o$($D$) $\geq$ 5 is odd, then span($V_{K}(t)$) $\leq$ n-3
\end{enumerate} \label{boundtheorem}
\end{theorem}
\begin{proof}
Suppose $D$ and $K$ are as given and let $G$ be the nonalternating skeleton for $D$.  For $o(D) = 3$ the only possible skeleton is a hexagon (from the required structure of $G$).  The extreme states bound gives the desired inequality.  The second case is the Kauffman-Murasugi result \cite{16}.  The third case follows from the previous lemma.
\end{proof}

Is this bound on the span sharp?  In one sense it is, as there are infinitely many prime knots for which the bound equals the span, as described below.  However, one can find prime knots where the weakness of the bound (that is, the difference $bound$ - $span$) can be made arbitrarily large.  For numerous examples of these, see \cite{9}.

Using the notation of \cite{12}, we call a diagram \textit{$+$ adequate} (in Thistlethwaite's sense) (resp. $-$ adequate) when $max(S_{A}) > max(S)$ (resp. $min(S_{B}) < min(S)$) for all states $S \neq S_{A}$ (resp. $\neq S_{B}$).

\begin{theorem}
\cite{14} If a diagram $D$ for a knot $K$ is $+$ adequate, then the extreme coefficient of maximal degree of $V_{K}(t)$ is $\pm1$.  If the diagram is $-$ adequate, then the extreme coefficient of minimal degree is $\pm1$. \label{adequatetheorem}
\end{theorem}

\begin{figure}[hbtp]
  \begin{center}
    \leavevmode
    \epsfxsize = 3.8cm
    \epsfysize = 3cm
    \epsfbox{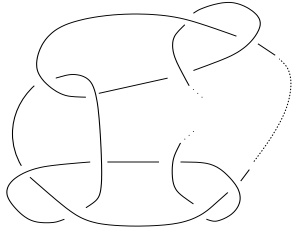}
    \caption{Knot $K_{n}$}
    \label{K11n74n}
  \end{center}
\end{figure}

Using the knot diagram $K = K11n74$, from the Hoste-Thistlethwaite table, construct the prime, reduced diagram $K_{n}$ shown in Fig. \ref{K11n74n} (the leftmost dotted lines loop 2$n$ additional times around the right dotted line, alternating as the lower two loops do).  One can check that the span of $V_{K_{n}}$($t$) equals the bound gotten from above as follows.  A simple inductive argument shows that for all $n$, $K_n$ is $+$ and $-$ adequate.  By Theorem \ref{adequatetheorem} the extreme coefficients are $\pm1$ so that the span equals the computed bound.  Furthermore, the nonalternating skeleton of $K_n$ has the following form: two $2(n+2)$-valent vertices connected by segments with one one interior vertex, yielding a total of $2(n+2) + 2$ vertices.  Therefore the bound is $c(D) - 1$, the maximal span for a nonalternating prime, reduced diagram with $c(D)$ crossings.

Knots of large enough crossing index or links with disjoint components could have disconnected skeletons.  It is natural to ask how this affects the Jones polynomials of such objects.  Let $\mid$$G$$\mid$ denote the number of connected components of $G$.

\begin{theorem}
Let $c(D) = n$, $o(D)=m$ and $\mid$G$\mid$ = p.  If $m$ is even, then span($V_{K}$($t$)) $\leq n-p$ and span($V_{K}$($t$)) $\leq n-p-1$ if $m$ is odd.\label{disconnect}
\end{theorem}
\begin{proof}
Suppose $G$ has components $G_1$, ..., $G_p$.  Maximizing the quantity \linebreak$q$ = $\mid$$G_A$$\mid$ + $\mid$$G_B$$\mid$ - $v$($G$) will result from maximizing $q_i$ on $G_i$.

$G_i$ will have $2m_i$ edges, with $m_{i} = 3$, $m_{i}$ odd and $m_{i} \geq 5$, or $m_i$ even.  In these cases, $q_{i} = -4$, $-6$, or $-2$, respectively.  If $m$ is even, then we may have all $q_i = -2$.  If $m$ is odd, we cannot have all $m_{i}$ even.  $q$ will obtain its maximum when $m_{i} = 3$ for only one $i$ and the remaining $m_{j}$ are even.  Hence, $q = q_{1}+...+q_{p} \leq -2p$ if $m$ is even and $q = \leq -2(p - 1)-4$ if $m$ is odd.  The result follows from Theorem \ref{extremestatesbound}.\\
\end{proof}

Notice for $m$ odd, regardless of $G$ having one or two components, \linebreak span$(V_{K}(t))$ $\leq n-3$.

\section{Classes of knots}
\label{classes}

From Theorem \ref{boundtheorem} we can conclude various facts about the Jones polynomials of certain classes of knots.  A \textit{satellite to a knot $K$} is gotten by embedding a knot in the solid $T$ and then tying the torus in the knot type $K$.  $K$ is referred to as a \textit{companion} to its satellite.  When tied, the $T$ may be twisted.  We can avoid this, however, by insisting that the longitude of $T$ be identified with the longitude of the companion knot.  We shall refer to such a knot as an \textit{untwisted satellite knot}.  

\begin{figure}[hbtp]
  \begin{center}
    \leavevmode
    \epsfxsize = 3cm
    \epsfysize = 2cm
    \epsfbox{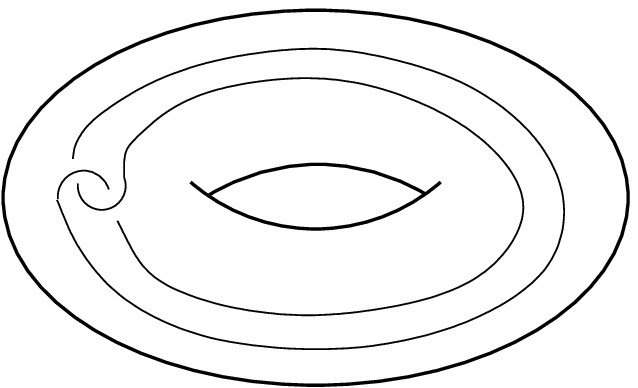}
    \caption{Constructing Whitehead doubles}
    \label{whitehead}
  \end{center}
\end{figure}

\begin{figure}[hbtp]
  \begin{center}
    \leavevmode
    \epsfxsize = 3cm
    \epsfysize = 2cm
    \epsfbox{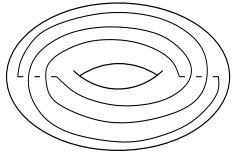}
    \caption{Creating $3_{2}$-strand cable knots}
    \label{3strandcable}
  \end{center}
\end{figure}

We consider two specific types of satellite knots: Whitehead doubles and $n_{m}$-cable knots.  They are formed via the unknotted circles lying in the tori as pictured in Figs. \ref{whitehead} and \ref{3strandcable}, respectively.  We will refer to the crossing regions of the unknotted circles lying in the torus (as well as their images upon forming the satellite knots) as the \textit{clasps of $K$}.

\subsection{Whitehead doubles}

\begin{theorem}
If $D$ is a totally reduced, nontrivial knot diagram for a knot $K$ and $K'$ is the untwisted Whitehead double for $K$, then span$(V_{K'}(t)) \leq 4c(D) - 1.$ \label{Whiteheadspan}
\end{theorem}
\begin{proof}
$D'$, the usual diagram for the standard Whitehead double of $K$, is totally reduced and has $c(D') = 4c(D) + 2$.  Moreover, $D'$ has four nonalternating edges (two positive and two negative) occurring at each region corresponding to a crossing in $D$.  We can assume the clasp of $D'$ is placed between two crossing regions that correspond to a nonpositive, nonnegative edge of $D$.  Therefore there are two nonalternating edges (one of each type) occuring at the clasp of $D'$.  Each positive or negative edge in $D$ corresponds to exactly two edges of the same type in $D'$. Hence, $o(D')$ is odd (specifically $o(D') = 2c(D) + 2o(D)+ 1$) and is greater than five (since $K$ is nontrivial).  Theorem \ref{boundtheorem} proves the result. 
\end{proof}

The above bound is exact for the untwisted Whitehead double of the trefoil, but for other alternating knots it appears to grow weaker as the crossing number of the original knot increases.  This is explained by the extreme states bound.  For the next few results, assume $D'$ is a untwisted Whitehead double for a nontrivial, totally reduced, alternating diagram $D$ of knot type $K$.

\begin{lemma}
$|S_{A}(D')| = 2|S_{A}(D)| - 1$ and $|S_{B}(D')| = 2|S_{B}(D)| + 1$.  \label{whiteheadstates}
\end{lemma}
\begin{proof}
Since $D$ is alternating, the regions of the complement of the shadow of $D$ can be 2-colored such that no two adjacent regions have the same coloring \cite{16}.  The state circles gotten from splitting $S_{A}$ and $S_{B}$ correspond precisely to the shaded and unshaded regions, respectively.  Without loss of generality, the shaded regions of $D$ have the form of a polygon with stacked edges, as in Fig. \ref{alternatingregions}(a).  These regions correspond to the portions of $D'$, as shown in Fig. \ref{alternatingregions}(b) and \ref{alternatingregions}(c) .  The clasp of $D'$ occurs in exactly one of these regions.

\begin{figure}[hbtp]
  \begin{center}
    \leavevmode
    \epsfxsize = 5cm
    \epsfysize = 4cm
    \put(58,58){(a)}
    \put(25,-15){(b)}
    \put(110,-15){(c)}
    \epsfbox{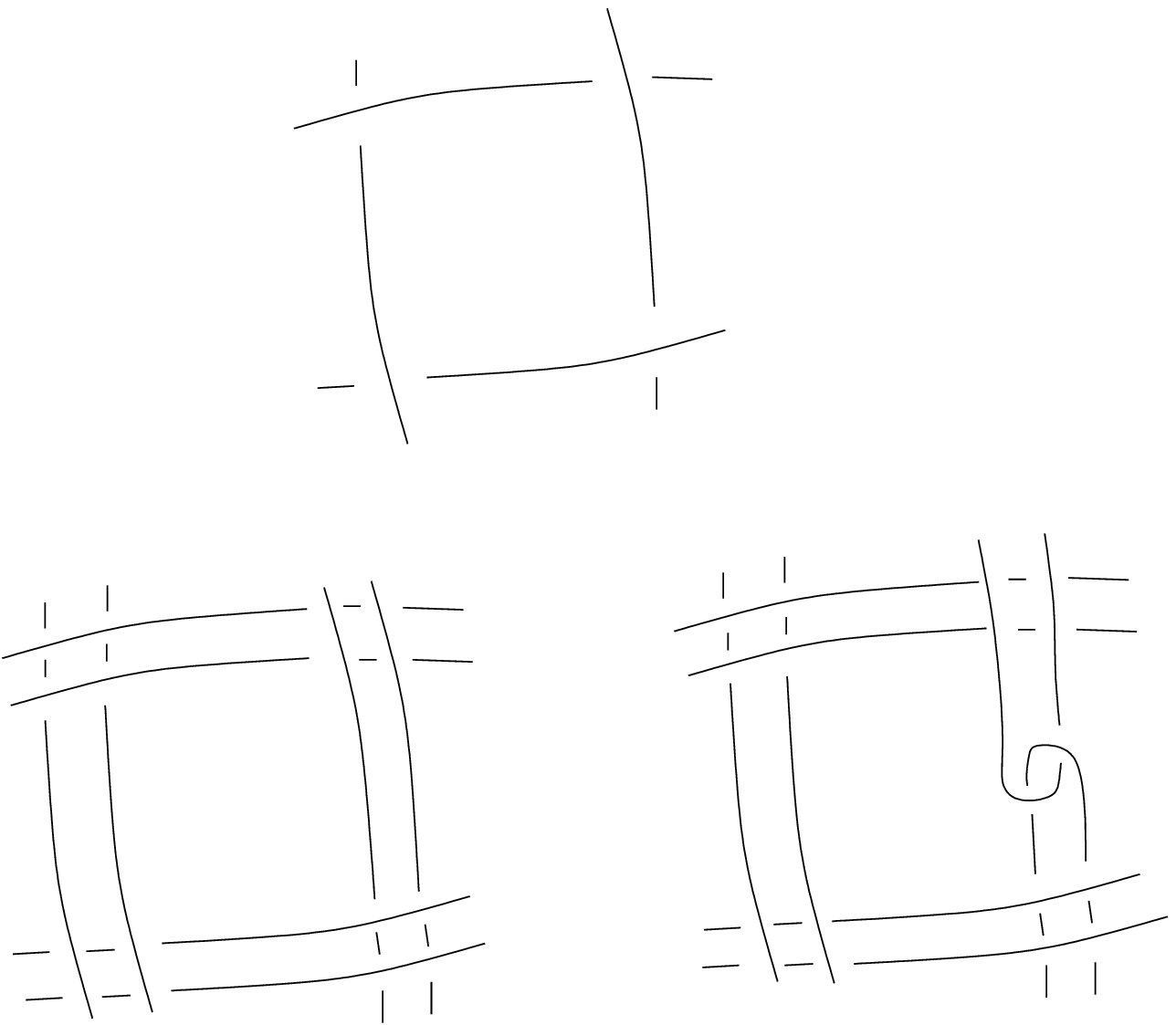}
    \caption{Alternating regions and their Whitehead doubles}
    \label{alternatingregions}
  \end{center}
\end{figure}

Consider $S_{A}(D')$.  Each region not containing the clasp is a specific case of a result found in \cite{8} and hence splits into exactly two state circles.  The region containing the clasp splits into exactly one state circle.  An argument involving $S_{B}(D')$ and unshaded regions follows similarly.  In this case, however, the clasp admits one extra circle while each region corresponds to exactly two state circles of $S_{B}(D')$.  Hence, the result follows.
\end{proof}

Before the main results on Whitehead doubles, we prove some preliminaries about the extreme coefficients of the Jones polynomial.  Recall that a set of vertices in a graph is called \textit{independent} if there are no edges connecting pairs of vertices from that set.

Let $f$ be a function defined on graphs $\Gamma$ as follows: $f(\Gamma) = \sum_{C}(-1)^{|C|}$, where the sum is taken over all independent sets of vertices of $\Gamma$.  Having formed the Lando graphs $L_{A}$ and $L_{B}$, we can calculate $a_{S_{A}}$ and $b_{S_{B}}$: \linebreak$a_{S_{A}} = (-1)^{|S_{A}|-1}f(L_A)$ ($b_{S_B}$ calculated similarly) \cite{2}.  Then $a_{S_A} \neq 0$ if and only if $f(L_A) \neq 0$.  The function $f$ has two key properties:

\begin{enumerate}
\item (Recursion) $f(\Gamma)$ =  $f(\Gamma - v)$ - $f(\Gamma - N(v))$, for a specific vertex $v$, where $\Gamma - v$ is the graph obtained from $\Gamma$ by deleting $v$ and all edges to which it is an endpoint and $\Gamma - N(v)$ results from deleting all neighboring vertices to $v$ (including $v$) and their adjoining edges.
\item (Disjoint Union) $f(\Gamma \coprod \Lambda)$ = $f(\Gamma) f(\Lambda)$. 
\end{enumerate}

\begin{lemma}
If $\Gamma$ is a complete $n$-partite graph, then $f(\Gamma) \neq 0$.  \label{npartite}
\end{lemma} 
\begin{proof}
Proceed by induction on $n$.  For $n = 2$, suppose $|v_i| = m_i$, where the $v_i$ are the vertex sets of the $n$-partite graph.  If $m_1$ = 1, call $v$ the one element of $v_1$. Then:
\begin{align*}
f(\Gamma) &= f(\Gamma - v) - f(\Gamma - N(v)) \\
&= 0 - 1 \\
&\neq 0 
\end{align*}

The first line of the above equality follows from the recursion propety.  The edge set of $\Gamma - v$ is empty and $\Gamma - N(v)$ is empty, giving the second line.  Now if $m_1 > 1$, choose $v$ to be an element of $v_1$.  Then as above, we have:

\begin{align*}
f(\Gamma) &= f(\Gamma - v) - f(\Gamma - N(v)) \\
&= f(\Gamma - v) - 0 \\
&\neq 0 
\end{align*}

$\Gamma - N(v)$ is a nonempty edgeless graph, proving the second equality.  By induction on $m_1$, the final equality follows.  Thus the result holds for $n = 2$.  For a general $n > 2$, inducting again on $m_1$ proves the result.
\end{proof}

Considering $D$, $D'$, $K$, and $K'$ as before, we have:

\begin{theorem}
span($V_{K'}(t)$) $\leq 3c(D) + 2$. \label{whiteheadjonesbound}
\end{theorem}
\begin{proof}
\begin{align*}
span\kb{D'} &\leq 2c(D') + 2(|S_{A}(D')| + |S_{B}(D')| - 2)\\
&= 8c(D) + 2(S_{A}(D') + S_{B}(D')) \\
&= 8c(D) + 2(2|S_{A}(D)| - 1 + 2|S_{B}(D)| + 1)\\
&= 8c(D) + 4(|S_{A}(D)| + |S_{B}(D)|)\\
&= 8c(D) + 4(c(D) + 2)\\
&= 12c(D) + 8
\end{align*}

The inequality is the extreme states bound (for $\kb{D}$).  The second equality follows from Lemma \ref{whiteheadstates} and the fourth equality is from \cite{16}.

Therefore span$\kb{D'}$ $\leq 12c(D) + 8$, implying that span($V_{K'}(t)$) $\leq 3c(D) \nolinebreak+\nolinebreak 2.$
\end{proof}

Notice that the nonalternating skeleton $G$ of $D'$ is isotopic to the shadow of $D$.

\begin{theorem}
$V_{K'}(t)$ has nonzero extreme coefficients and hence span equal to $3c(D) + 2$.  \label{whiteheadcoefficients}
\end{theorem}
\begin{proof}
Except for near the clasp of $D'$, $G$ lies entirely in the interior of $D$.  Figure \ref{nonclaspwhiteheadlandoarcs} shows how arcs behave around a section of $D'$ corresponding to a crossing of $D$ but not near the clasp of $D'$.  Notice that both endpoints of no arc lie on the same component of $G_{A}$ or $G_{B}$ (recall Fig. \ref{alternatingregions}($b$).

\begin{figure}[hbtp]
  \begin{center}
    \leavevmode
    \epsfxsize = 4cm
    \epsfysize = 5cm
    \put(72,125){$A$}
    \put(72,87){$B$}
    \put(18,60){$G_{A}$}
    \put(80,60){$G_{B}$}
    \epsfbox{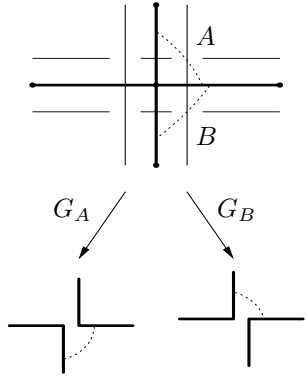}
    \caption{Arcs not near the clasp of $D'$}
    \label{nonclaspwhiteheadlandoarcs}
  \end{center}
\end{figure}

Now consider a neighborhood of the clasp that includes the crossings around it.  We need only consider the arcs whose endpoints lie on the portion of $G$ that exits the interior of $D'$.  Two of these arcs ($e_{1}$ and $e_{2}$) come from the clasp crossings and one from each of the sections of $D'$ corresponding to crossings in $D$, and all are of the same type.  See Fig. \ref{whiteheadclasparcs}.

\begin{figure}[hbtp]
  \begin{center}
    \leavevmode
    \epsfxsize = 4.5cm
    \epsfysize = 4.5cm
    \epsfbox{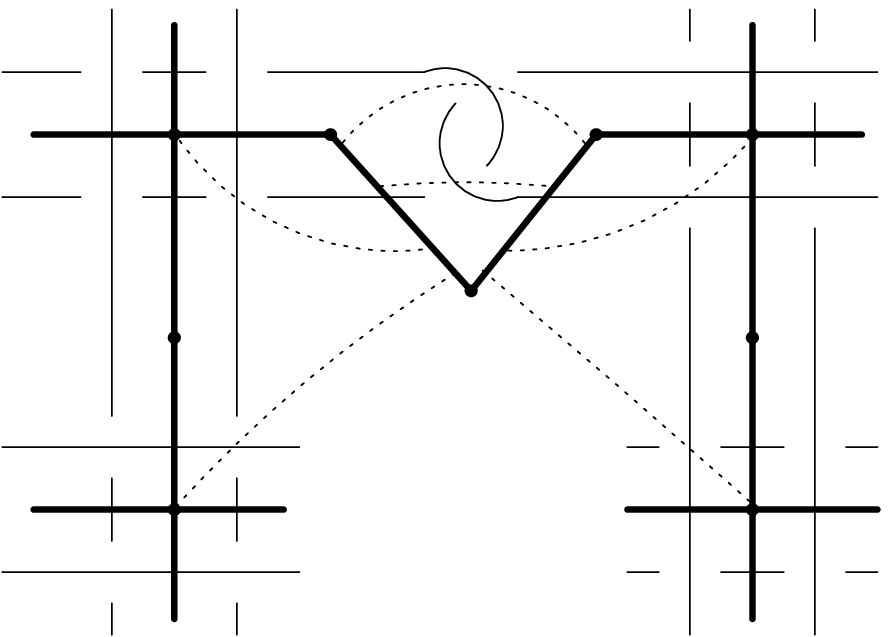}
    \caption{Near the clasp of $D'$}
    \label{whiteheadclasparcs}
  \end{center}
\end{figure}

Therefore one of the Lando graphs is bipartite and the other is empty.  By Lemma \ref{npartite}, then, both extreme coefficients of $V_{K'}(t)$ are nonzero. 
\end{proof}

Theorem \ref{whiteheadcoefficients} partially answers an unsolved question given in \cite{1}:

\begin{corollary}
If $K'$ is an untwisted Whitehead double to a nontrivial, alternating knot $K$, then $c(K) < c(K')$.
\end{corollary}
\begin{proof}
Take $D$ to be a totally reduced, alternating diagram for $K$ and let $D'$ be an untwisted Whitehead double diagram to $D$.  Then, 
\begin{align*}
span(V_{K}(t)) &= c(K)\\
&= c(D)\\
&< 3c(D) + 2\\
&= span(V_{K'}(t))\\
&\leq c(K')
\end{align*}
\end{proof}

The proof of Theorem \ref{whiteheadcoefficients} provides even more results for any Whitehead double of an alternating knot.  

\begin{corollary}
Any Whitehead double of an alternating knot has nontrivial Jones polynomial and hence is a nontrivial knot. \label{alternatingwhitehead}
\end{corollary}
\begin{proof}
Take $D$ to be a totally reduced, alternating, nontrivial knot diagram for a knot $K$ and $K'$ the knot type of any Whitehead double $D'$ for $D$.  An arc of $D$ corresponds to a pair of parallel arcs in $D'$.  A general diagram for $D'$ may have twists, either positive or negative, between any such pair.  Notice, however, that some of these twists may be removed (by simply untwisting the edges in $K'$, leaving an equivalent diagram with only one type of twist.  Moreover, we can slide (by stretching the knot edges) these twists so that they all occur between two sets of crossings of $D'$ corresponding to adjacent crossings in $D$.  Finally, we can assure ourselves that these parallel edges containing the twists do not share a complementary region of the plane with the clasp.

If twists occur between parallel arcs in $D'$, then $o(D') = o(D) + 1$, with the extra two edges of the nonalternating skeleton occurring near the twists.  As in the clasp case of the proof of Theorem \ref{whiteheadcoefficients}, all Lando arcs must be of the same type.  Therefore at most one of the extreme coefficients of $V_{K'}(t)$ differs from the nonzero extreme coefficients of $V_{K}(t)$.  Thus, $V_{K'}(t)$ is nontrivial.
\end{proof}

\subsection{Cable knots}
Most of the results for Whitehead doubles carry over similarly to untwisted $n_{1}$-strand cable satellites.  Assume $D'$ is a untwisted $n_{1}$-strand cable diagram for a totally reduced, nontrivial, alternating diagram $D$ of a knot $K$.  For the results and proofs below, the roles of $A$ and $B$ may be reversed, depending on whether the clasp connects via an arc over or under the other $n - 1$ arcs.  

\begin{figure}[hbtp]
  \begin{center}
    \leavevmode
    \epsfxsize = 5cm
    \epsfysize = 6cm
    \put(48,72){$S_{A}$}
    \epsfbox{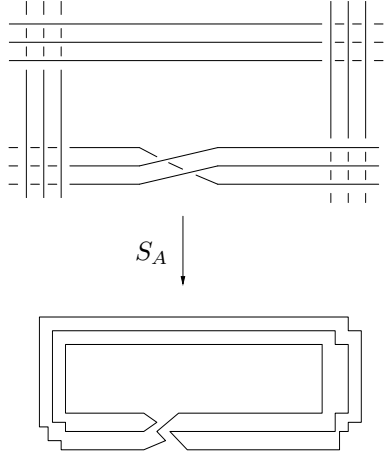}
    \caption{$S_{A}$ for $3$-clasp knot}
    \label{nclaspstate}
  \end{center}
\end{figure}

\begin{lemma}
$|S_{A}(D')| = n|S_{A}(D)|$ and $|S_{B}(D')| = n(|S_{B}(D)| - 1) + 1$. \label{nstrandstatecount}
\end{lemma}
\begin{proof}
As we did with Whitehead doubles, we consider the two types of regions in $D'$: those with and those without the clasp on their boundary.  The states circles of $S_{A}(D')$ and $S_{B}(D')$ correspond to the state circles of $S_{A}(D)$ and $S_{B}(D)$, respectively.  The regions of $D'$ without the clasp on their boundary result in exactly $n$ state circles for each corresponding closed curve in $S_{A}(D)$, as do those in $S_{B}(D')$.  Exactly one state circle of $S_{A}(D)$ will correspond to the clasp of $D'$ and exactly one corresponding circle results when forming $S_{A}(D')$.  See Fig. \ref{nclaspstate}  A similar result holds for $S_{B}(D')$, with $n$ state circles being formed around the clasp.
\end{proof}

\begin{theorem}
$V_{K'}(t)$ has nonzero extreme coefficients and hence span equal to $\frac{1}{2}(n^{2} + n)c(D) + (n-1)$.
\end{theorem}
\begin{proof}
The proof follows exactly as the proof to Theorem \ref{whiteheadcoefficients}.  Notice that $c(D') = n^{2}c(D) + n - 1$.  
\end{proof}

As with Whitehead doubles, we have proven a case of Adams' unsolved question:

\begin{corollary}
If $K'$ is an untwisted $n_{1}$-strand cable satellite to a nontrivial, alternating knot $K$, then $c(K) < c(K')$.
\end{corollary}

Finally, we consider the standard $n_{m}$-strand cable satellite diagram $D'$ of an alternating knot for any $m$.  Notice that we can stretch the edges of $D'$ so that all clasps are adjacent (just as in the proof of Corollary \ref{alternatingwhitehead}).  Also note that twists in the torus can be realized by adding or removing clasps (since we only allow knots to be formed via twists, $m$ will remain relatively prime to $n$).  Thus it is enough to consider only untwisted $n_{m}$-strand cable knots.  Just as with the Whitehead doubles to alternating knots, only one type of Lando chord will be preserved after splitting into $G_{A}$ and $G_{B}$.  Therefore either $L_{B}$ or $L_{A}$ is empty, yielding at least one nonzero extreme coefficient of $V_{K'}(t)$, proving:

\begin{theorem}
Any $n_{m}$-strand cable knot of an alternating knot has nontrivial Jones polynomial and hence is a nontrivial knot.
\end{theorem}

\subsection{Polygonal skeletons}
Another class of knot diagrams can be classified by their skeletons.  Consider the case where $G$ is a 2$p$-gon skeleton, $p \geq 2$, for a prime and totally reduced diagram.  Then $D$ has 2$p$ non-alternating edges.  These are the only edges of $D$ that intersect $G$ when $G$ is superimposed on the diagram $D$. For a given non-alternating edge $e$, call the two edges of $D$ which intersect adjacent edges of $G$ the \textit{neighbors} of $e$.

\begin{figure}[hbtp]
  \begin{center}
    \leavevmode
    \epsfxsize = 3cm
    \epsfysize = 2cm
    \epsfbox{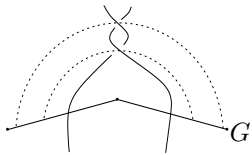}
    \put(0,5){$G$}
    \caption{Neighbors contributing Lando-$A$ arcs}
    \label{contributinglandoarcs}
  \end{center}
\end{figure}

The vertex conditions on $G$ prove that the neighbors of \textit{e} necessarily have the same sign.  Suppose an edge crosses the superimposed skeleton.  If the crossing it first encounters is a twist with one of its neighbors, then a Lando arc must exist at thost twists.  See Fig. \ref{contributinglandoarcs}.  $n$ half twists between neighbors results in $n$ Lando arcs of the same type.

If an edge $e$ and its neighbor $e_{n}$ form a twist as just described and the resulting arc is a Lando-$A$ arc (that is, passing through $A$ regions), we will say that \textit{$e$ contributes a Lando-A arc with $e_{n}$ to $L_{B}$}.  An edge and its neighbor similarly contribute Lando-$B$ arcs to $L_{A}$.   There are restrictions, however, on how pairs of neighbors interact, given by the following lemmas.  Recall that $G$ is considered to be a polygon and hence splits the plane into two distinct regions.

\begin{lemma}
Suppose an edge \textit{e} and one of its neighbors \textit{$e_n$} contribute a Lando arc to $L_A$ ($L_B$).  An adjacent pair of neighbors then cannot contribute Lando arcs to $L_B$ ($L_A$) in the same region that \textit{e} and \textit{$e_n$} contribute their Lando arc. \label{adjacentlandosame}
\end{lemma}
\begin{proof}
Suppose \textit{e} and one of its neighbors \textit{$e_n$} contribute a Lando arc to $L_A$.  The other neighbor to \textit{e} must twist in an opposite order with its other neighbor in order to contribute a Lando arc to $L_B$.  This contradicts the vertex structure of $G$. 
\end{proof}

Similarly:

\begin{lemma}
Suppose an edge \textit{e} and one of its neighbors \textit{$e_n$} contribute a Lando arc to $L_A$ ($L_B$).  An adjacent pair of neighbors then cannot contribute Lando arcs to $L_A$ ($L_B$) in the opposite region that \textit{e} and \textit{$e_n$} contribute their Lando arc. Moreover, \textit{e} also contributes to a Lando arc with its other neighbor only if they contribute it to $L_A$ in the opposite region in which \textit{e} and $e_n$ contribute. \label{adjacentlandodifferent}
\end{lemma}

\begin{figure}[hbtp]
  \begin{center}
    \leavevmode
    \epsfxsize = 4cm
    \epsfysize = 2cm
    \epsfbox{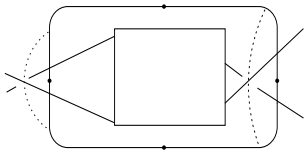}
    \caption{Removable Lando arcs}
    \label{removablelandoarcs}
  \end{center}
\end{figure}

These two lemmas lead to the following theorem:

\begin{theorem}
If $D$ is a totally reduced diagram for a knot $K$ and $o(D) = 2$, then $V_{K}(t)$ is nontrivial.  In particular, the result holds for any $K$ with $o(K) = 2$.  \label{olength2}
\end{theorem}
\begin{proof}
If $o$($K$) = 2, then $K$ has a totally reduced diagram $D$ of minimal $o$-length two.  Take $G$ to be the non-alternating skeleton of $D$.  $G$ is a quadrilateral.  If there is a Lando-$A$ arc $a$ in one of the regions formed by $G$, then there are no Lando-$B$ arcs in that same region.  By Lemma \ref{adjacentlandodifferent}, any Lando-$B$ arcs must be formed by the nonalternating edges not used to form $a$. These are removable through a twist; see Fig. \ref{removablelandoarcs}.  Thus we can find a diagram for $D$ where only one type of Lando chord is present.  Therefore one of the Lando graphs for $D$ is empty, proving that $V_{K}(t)$ is nontrivial.
\end{proof}

Recall that an \textit{almost alternating diagram} is a nonalternating knot diagram that becomes an alternating diagram upon switching the over- and undercrossing strands at exactly one crossing.  An \textit{almost alternating knot} is a  knot that has an almost alternating diagram.  The following corollary follows immediately.

\begin{corollary}
Every almost alternating knot has nontrivial Jones polynomial.
\end{corollary}

Lemma \ref{npartite} proves that any diagram with an $n$-partite Lando graph $L_A$ or $L_B$ must have a nonzero extreme coefficient and hence nontrivial Jones polynomial.  We classify some of these in the following theorem:

\begin{theorem}
Suppose $D$ is a diagram for a knot $K$ such that the nonalternating skeleton $G$ for $D$ is a single polygon.  If the only Lando arcs of $D$ arise from twists between neighbors where any arc contributed by neighbors $e$ and $e_n$ has one of $e$ or $e_n$ contributing an arc with its other neighbor, then $V_{K}(t)$ is nontrivial. \label{npartiteknots}
\end{theorem}
\begin{proof}
The Lando graphs for such a diagram are either disjoint unions of complete $n$-partite graphs or a complete $n$-partite graph such that \textit{all} the disjoint sets of vertices are completely connected (including $v_1$ completely connected to $v_n$).  The first case follows from Lemma \ref{npartite} and the disjoint union property of $f$.  Since $G$ is a single polygon, the latter case will occur only if $D$ has no Lando arcs of one type.  Therefore in this case, one of the Lando graphs is empty and the result follows.
\end{proof}

\begin{figure}[hbtp]
  \begin{center}
    \leavevmode
    \epsfxsize = 8cm
    \epsfysize = 4cm
    \put(190,73){$T_{1}$}
    \put(117,40){$T_{2}$}
    \put(61,90){$T_{3}$}
    \put(25,21){$T_{4}$}
    \epsfbox{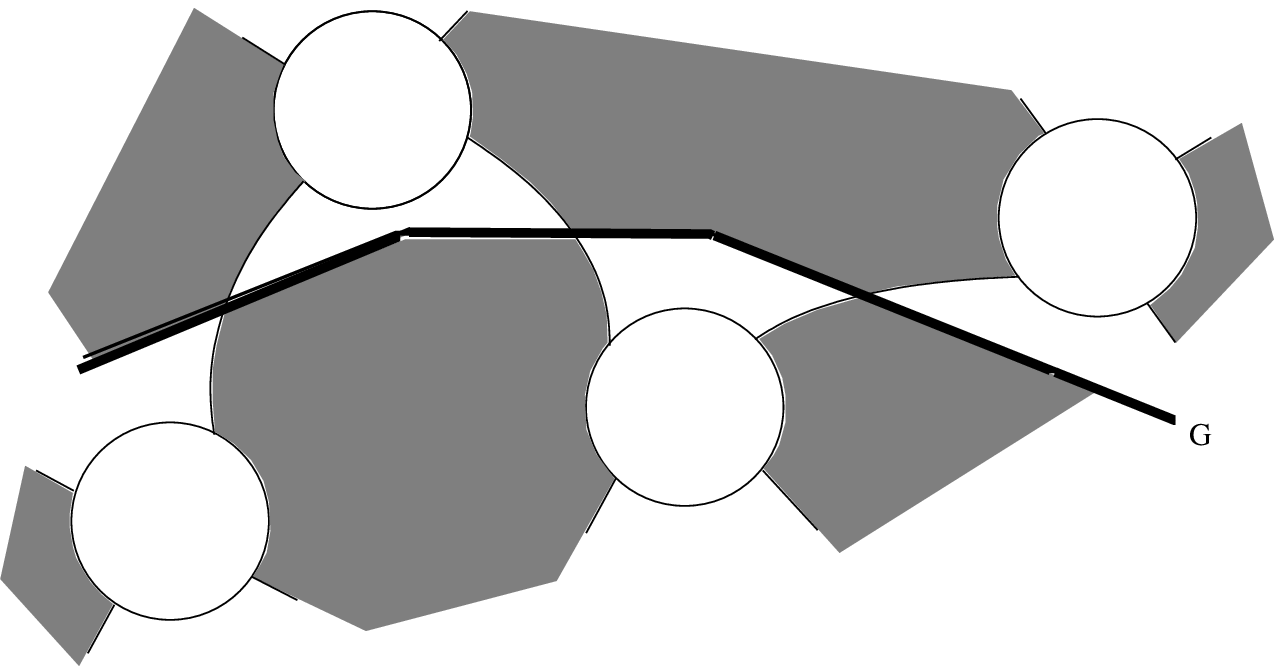}
    \caption{Illustrating Theorem \ref{npartiteknots}}
    \label{polytwists}
  \end{center}
\end{figure}

\begin{figure}[hbtp]
  \begin{center}
    \leavevmode
    \epsfxsize = 4cm
    \epsfysize = 4cm
    \epsfbox{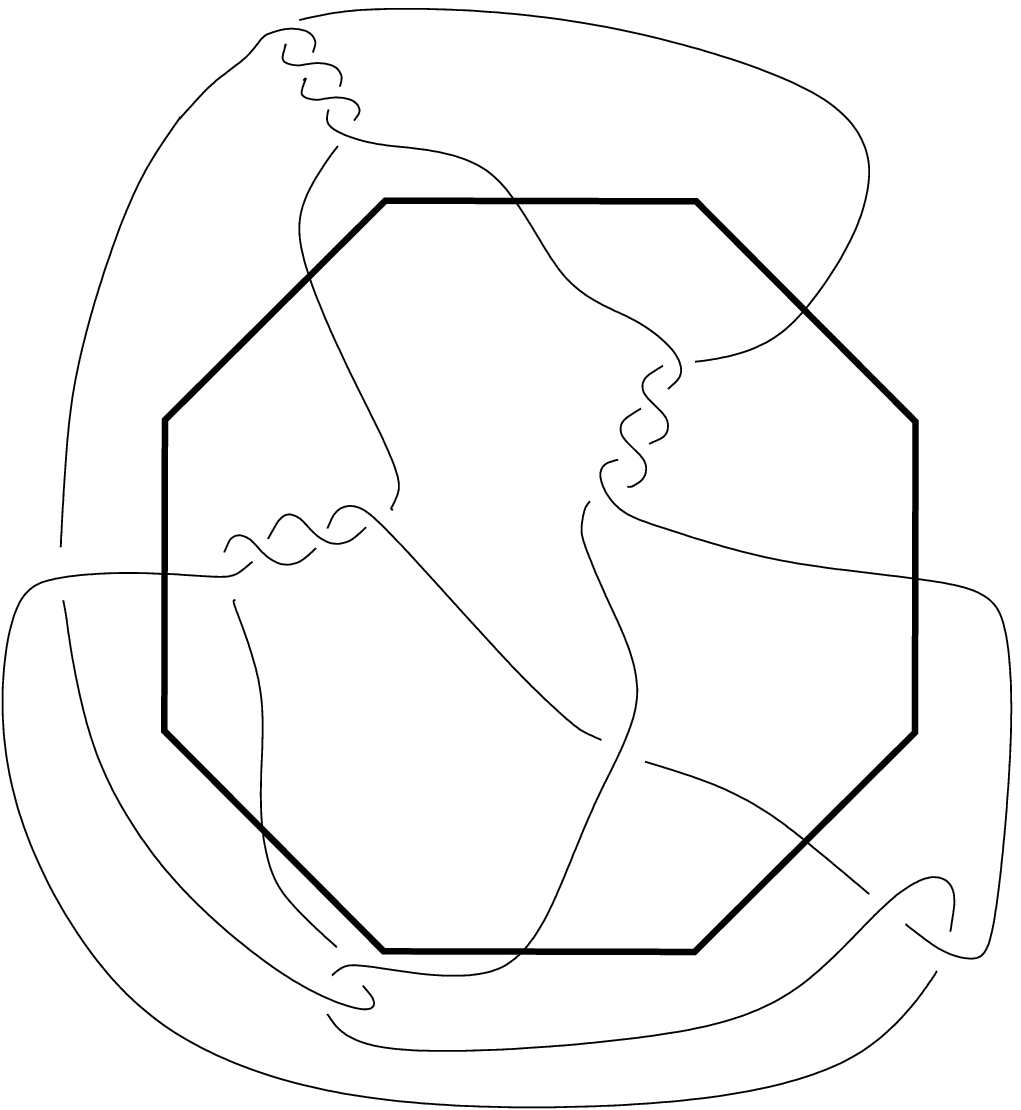}
    \caption{Example of Thereom \ref{npartiteknots}}
    \label{twistexample}
  \end{center}
\end{figure}

Figure \ref{polytwists} illustrates the knots to which Theorem \ref{npartiteknots} applies.  The 2-tangles consist of allowable twists (see Lemmas \ref{adjacentlandodifferent} and \ref{adjacentlandosame}), and no crossing outside these tangles bounds shaded regions across it. Figure \ref{twistexample} is an example of such a knot.

\subsection{Pretzel knots}
We consider one last class of knots: pretzel knots.

\begin{figure}[hbtp]
  \begin{center}
    \leavevmode
    \epsfxsize = 7cm
    \epsfysize = 3cm
    \put(14,39){$c_{1}$}
    \put(72,39){$c_{2}$}
    \put(175,39){$c_{n}$}
    \epsfbox{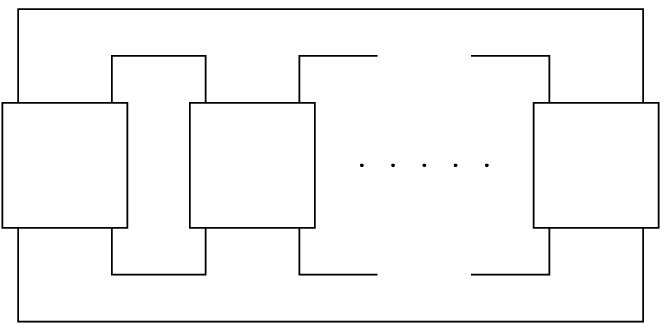}
    \caption{$P(c_{1}, c_{2}, ..., c_{n})$}
    \label{pretzel}
  \end{center}
\end{figure}

A \textit{pretzel link} $P(c_{1}, c_{2}, ..., c_{n})$ is a link of the type pictured in Fig. \ref{pretzel}, where the boxes represent vertically twisted $2$-tangles.  If $c_{i} < 0$ then the twists are negative (the slope of the overcrossing strand is negative in the usual planar sense); if $c_{i} > 0$ then the tangle consists of $c_{i}$ positive half twists.   It is a link of at most two components.  Moreover, it is a knot if and only if $n$ is odd and $c_{i}$ are odd for all $i$, or $c_{i}$ is even for just one $i$ \cite{7}.  We will assume from this point that all mentioned pretzel links are indeed knots.

\begin{lemma}
If $P(c_{1}, c_{2}, ..., c_{n})$ is a pretzel knot with $D$ its usual diagram, then $o(D)$ equals the number of changes in sign from $c_{i}$ to $c_{i+1}$ (considered cyclically).
\end{lemma}
\begin{proof}
Like rational knots, the only nonalternation in pretzel knots can occur between twists of opposite signs.  Adjacent twists of opposite sign contribute one positive edge and one negative edge to $D$.
\end{proof}

We know sets of distinct knots with the same Jones polynomial ($P(a_{1}, a_{2}, ..., a_{n})$ is isotopic to $P(a_{j+1}, ..., a_{n}, a_{1}, ..., a_{j})$ but not necessarily to any permutation of the $a_{i}$, and the Jones polynomial of $P$ is not changed over any permutation of the $a_{i}$).  But can we conclude anything about the triviality of these polynomials?  In \cite{8}, Lickorish and Thistlethwaite show that $P(a_{1}, a_{2}, ..., a_{n}, b_{1}, b_{2}, ..., b_{m})$ has nontrivial Jones polynomial, where $a_{i} \geq 2$, $-b_{j} \geq 2$, $n \geq 2$, and $m \geq 2$.  Landvoy exhibits an infinite family of pretzel knots, any of the form \linebreak$P(c_{1}, -2c_{1}-1, -2c_{1}+1)$, for which their Jones polynomials are nontrivial \cite{7} yet, from Parris, all have trivial Alexander polynomial \cite{10}.  Landvoy's result on the Jones polynomial follows immediately from Corollary \ref{olength2}. The nonalternating skeleton provides the following even stronger result.

\begin{theorem}
Any pretzel knot $P(c_{1}, c_{2}, ..., c_{n})$ with $|c_{i}|$ = 1 for at most one $i$ has nontrivial Jones polynomial. \label{pretzelthm}
\end{theorem}
\begin{proof}
Suppose $K = P(c_{1}, c_{2}, ..., c_{n})$ is a pretzel knot.  Either $c_{i}$ and $c_{i+1}$ are of different sign for all $i$ (considered cyclically) or an adjacent pair of these integers have the same sign.

Suppose the first case is true.  Note that $n$ must be even, since the alternation $c_{n}$ and $c_{1}$ must be of opposite sign.

For $n = 2$, by our previous remarks we have four edges that are positive or negative, so that $o(K) = 2$.  By Theorem \ref{olength2}, $V_{K}(t)$ is nontrivial.

Now consider a general postive value of $n = 2m$, $m > 1$.  We have two subcases, whether or not $|c_{i}| = 1$ for some $i$.  First suppose $|c_{i}| > 1$ for all $i$.  A result by Lickorish and Thistlethwaite \cite{8} gives $V_{K}(t) = V_{K'}(t)$, where up to signs $K' = P(c_{1}, c_{3}, ..., c_{2m-1}, c_{2}, c_{4}, ..., c_{2m})$.  By Lickorish and Thistlethwaite's result, $V_{K}(t)$ is nontrivial.

A second way of proving this first subcase involves Lando arcs and will be used in our next subcase.  Consider $c_{i}$ and $c_{i+1}$, $1 < i < n$, as in Fig. \ref{adjacenttwists}.  

\begin{figure}[hbtp]
  \begin{center}
    \leavevmode
    \epsfxsize = 4cm
    \epsfysize = 3cm
    \epsfbox{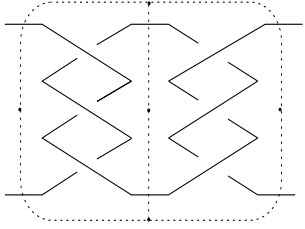}
    \caption{Adjacent twists with part of skeleton superimposed}
    \label{adjacenttwists}
  \end{center}
\end{figure}

Notice part of the nonalternating skeleton superimposed and how its edges inherit signs, as marked, from $K$.  Draw in the Lando arcs passing through the twists $c_{i}$ and $c_{i+1}$.  Lando-$B$ arcs pass through positive twists while negative twists contribute $A$ arcs.  Form the skeleton $G$ form $G_{A}$ and $G_{B}$.  Figure \ref{pretzellandoarcs} shows $G_{A}$ superimposed with the $B$-chords present.  Notice that their endpoints do not lie on the same component of $G_{A}$, giving no arcs contributed to $L_{A}$.  The same is true for $L_{B}$.

\begin{figure}[hbtp]
  \begin{center}
    \leavevmode
    \epsfxsize = 4cm
    \epsfysize = 3cm
    \epsfbox{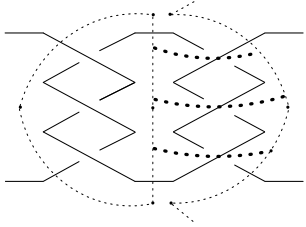}
    \caption{$G_{A}$ and Lando $B$-arcs}
    \label{pretzellandoarcs}
  \end{center}
\end{figure}

We must consider $c_{1}$ and $c_{n}$, however.  But the same result holds in this case.    Therefore $L_{A}$ and $L_{B}$ are empty, proving that $V_{K}(t)$ is nontrivial.

We move on to our second subcase, utilizing the Lando graphs.  As above, most possible Lando arcs have endpoints lying on different components of $G_{A}$ or $G_{B}$.  The only exception occurs when $|c_{i}| = 1$, as this twist has Lando arcs passing through both its $A$ and $B$ regions.  One of these arcs will have its endpoints lying on two distinct components of the graphs $G_{A}$ or $G_{B}$.  However, the Lando-$A$ arc passing through the $c_{i} = 1$ twist has its endpoints lying on a single component of $G_{B}$, proving that $L_{B}$ is nonempty.  A similar result holds with $L_{A}$ and the $-1$ twist.  This is the only case creating nonempty Lando graphs, and by our hypothesis, only one of these cases can hold.  Thus, either $L_{A}$ or $L_{B}$ is empty and consequently, $V_{K}(t)$ is nontrivial.

The first case has been proven.  If $c_{i}$ and $c_{i+1}$ have the same sign for some $i$, we may assume (by isotopy) that it is $c_{n}$ and $c_{1}$ with the same sign.  But again, nonalternation occurs only between twists of different signs and the proof follows as the first case, resulting in $L_{A}$ or $L_{B}$ (possibly both) empty.
\end{proof}

We have the following immediate corollary to Theorem \ref{pretzelthm}:

\begin{corollary}
Suppose $K = P(c_{1}, c_{2}, ..., c_{n})$ is a pretzel knot such that all $c_{i}$ with $|c_{i}|$ = 1 are of the same sign.  Then $V_{K}(t)$ is nontrivial. \label{samesignpretzel}
\end{corollary}
\begin{proof}
The proof follows the proof of Theorem \ref{pretzelthm}.  Only one of the Lando graphs is nonempty.
\end{proof}

\begin{figure}[hbtp]
  \begin{center}
    \leavevmode
    \epsfxsize = 6cm
    \epsfysize = 6cm
    \epsfbox{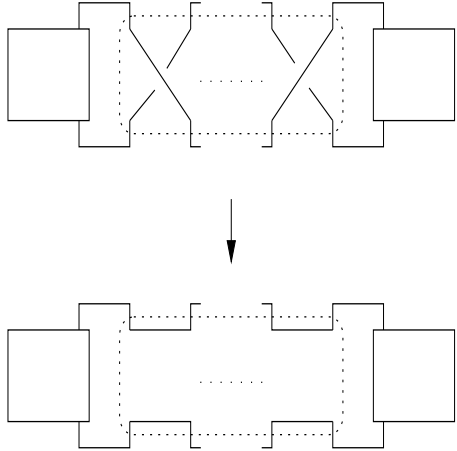}
    \caption{Removable twists}
    \label{pretzeltwist}
  \end{center}
\end{figure}

As a result of our above work, we have one of our main results:

\begin{theorem}
Every pretzel knot has nontrivial Jones polynomial. \label{everypretzel}
\end{theorem}
\begin{proof}
Let $K = P(c_{1}, c_{2}, ..., c_{n})$ be a pretzel knot.  It suffices to show it takes the form given in Corollary \ref{samesignpretzel}.  If $c_{i} = 1$ and $c_{j} = -1$ ($i < j$), then $K$ is isotopic to $P(c_{1}, c_{2}, ..., c_{i-1}, c_{i+1}, ..., c_{j-1}, c_{j+1}, ..., c_{n})$ by a simple twist, shown in Fig. \ref{pretzeltwist}.  Repeating yields a pretzel knot of the desired form.
\end{proof}

Theorem \ref{pretzelthm} shows for certain pretzel knots (namely those with $|c_{i}| > 1$ for all $i$) it is not difficult to calculate the actual span of the Jones polynomial.  For $K$ given in the first case of the proof of Theorem \ref{pretzelthm}, $o(K)$ is even.  Then we have span($V_{K}(t)$) = $c(D) - 1$ = $\sum |c_{i}| - 1$.  For those $K$ not of this type, count the number of sign changes between $c_{i}$ and $c_{i+1}$ (by isotopy we can assume $c_{1}$ and $c_{n}$ have the same sign).  This is the $o$-length $o$ for our diagram.  If $o = 0$, then the knot is alternating.  Else if $o$ is even, then span($V_{K}(t)$) = $\sum |c_{i}| - 1$, and if $o$ is odd, then span($V_{K}(t)$) = $c(D) - 1$ = $\sum |c_{i}| - 2$.

The proof of Theorem \ref{pretzelthm} also shows that $K$ (with $|c_{i}| > 1$) has extreme coefficients of $\pm$1 since $L_{A}$ and $L_{B}$ are empty.  Hence $a_{S_{A}} = (-1)^{|S_{A}|-1}$ and $b_{S_{B}} = (-1)^{|S_{B}|-1}$.  Consider how forming $S_{A}$ affects positive and negative twists.  Each positive $c_{i}$ contributes to two simple closed curves of $S_{A}$, which we count twice when considering all $i$.  A negative $c_{i}$ contributes $c_{i} - 1$ simple closed curves to $S_{A}$ and contribute to at most two more circles.  But if we assume all $c_{i}$'s are not of the same sign, then it will contribute exactly one additional circle.   Similar results hold for $S_{B}$.  If $p_{c}$ and $n_{c}$ are the number of positive and negative $c_{i}$, respectively, then $|S_{A}|$ = $p_{c}$ + $\sum -(c_{i}+1)$ (where we sum over all negative $c_{i}$) and $|S_{B}|$ = $n_{c}$ + $\sum (c_{i}-1)$ (summing over all $c_{i}$ positive).  In lieu of Theorem \ref{everypretzel}, the following lemma applies to all pretzel knots:

\begin{corollary} If $K = P(c_{1}, c_{2}, ..., c_{n})$ has $c_{i} \neq -1$ for all $i$, then \linebreak$a_{S_A} = (-1)^{p_{c} - 1 + \sum -(c_{i}+1)}$.  If $c_{i} \neq 1$ for all $i$, then $b_{S_B} = (-1)^{n_{c} - 1 + \sum (c_{i}-1)}$.
\end{corollary}

\section*{Acknowledgments}
The author is grateful to Raymond Lickorish for helpful comments provided on an earlier version of this paper.

\end{document}